\documentclass[12pt]{article}

\title{Finite Variation of Fractional L\'evy Processes}
\author{Christian Bender${}^1$ \and Alexander Lindner${}^2$ \and Markus Schicks${}^2$ }
\date{}

\usepackage{amssymb}
\usepackage{amsfonts,dsfont}
\usepackage{amsmath}
\usepackage{amsthm}
\usepackage{graphics}
\usepackage{latexsym}

\newtheorem{theorem}{Theorem}[section]
\newtheorem{lemma}[theorem]{Lemma}
\newtheorem{remark}[theorem]{Remark}

\newtheorem{proposition}[theorem]{Proposition}

\newcommand{\bR}{\mathbb{R}}
\newcommand{\bE}{\mathbb{E}}

\newcommand{\cI}{\mathcal{I}}
\newcommand{\cF}{\mathcal{F}}
\newcommand{\cB}{\mathcal{B}}

\newcommand{\One}{\mathds{1}}

\newcommand{\Var}{\mbox{\rm Var}\,}
\newcommand{\sign}{{\rm sign}}

\begin{document}
\footnotetext[1]{Universit\"at des Saarlandes, Fachrichtung
Mathematik, P.O Box 151150, 66041 Saarbr\"ucken, Germany. Email:
bender@math.uni-sb.de} \footnotetext[2]{Technische Universit\"at
Braunschweig, Institut f\"ur Mathematische Stochastik,
Pockelsstra{\ss}e 14, 38106 Braunschweig, Germany. Email:
\{a.lindner,\; m.schicks\}@tu-bs.de}
 \maketitle

\abstract{Various characterizations for fractional L\'evy process to
be of finite variation are obtained, one of which is in terms of the
characteristic triplet of the driving L\'evy process, while others
are in terms of differentiability properties of the sample paths. A
zero-one law and a formula for the expected total variation is also
given.}
\bigskip




\section{Introduction}

Recently there has been increased interest in fractional L\'evy
processes, which are generalizations of fractional Brownian motion.
Benassi et al.~\cite{Be1, Be} and Marquardt~\cite{Ma} introduced
real harmonizable fractional L\'evy processes, well-balanced (moving
average) fractional L\'evy processes $N_d$ and non-anticipative
(moving average) fractional L\'evy processes $M_d$. Apart from a
normalizing constant, these arise by replacing the Brownian motion
in the corresponding representation of fractional Brownian motion by
a centered square-integrable L\'evy process, and the precise
definitions of $M_d$ and $N_d$ are given below in \eqref{def-flp}
and \eqref{def-Nd}, respectively. Note that although the different
representations all give fractional Brownian motion if the driving
process is Brownian motion, in general the corresponding definitions
lead to different processes for arbitrary driving L\'evy processes.
However, all the mentioned processes have the same second order
structure as fractional Brownian motion. Other properties, such as
self-similarity, are not necessarily shared with fractional Brownian
motion (cf.~\cite{Be, Ma}), and in this paper we shall concentrate
on the semimartingale property and on finite variation of the paths.
 While it is well
known that a fractional Brownian motion with Hurst parameter $H\in
(0,1)$ cannot be a semimartingale unless $H=1/2$ (e.g.
Mishura~\cite{Mi}, p.~71), in particular cannot be of finite
variation on compacts almost surely for any $H\in (0,1)$, this is
not the case for $M_d$ and $N_d$. Marquardt~\cite{Ma}, Theorems~4.6,
4.7, has shown that $M_d$ will be of finite variation if the driving
L\'evy process is a compensated compound Poisson process, and has
given examples when $M_d$ is not a semimartingale. The aim of this
paper is to provide a complete characterization for the
non-anticipative fractional L\'evy process $M_d$ and for the
well-balanced fractional L\'evy process $N_d$ to be of finite
variation on compacts, equivalently for them to be semimartingales,
in the long memory case, i.e. when $H:= d+1/2\in (1/2,1)$. This
subject and the obtained results are closely related to a recent
paper by Basse and Pedersen~\cite{Ba}, who characterised the
semimartingale property of general one-sided L\'evy driven moving
average processes and applied their results to obtain a necessary
condition for the non-anticipative fractional L\'evy process $M_d$
to be of finite variation (see Remark~\ref{rem-Basse-Pedersen}
below). This condition is expressed in terms of the absence of a
Brownian motion component and an integrability condition on the
L\'evy measure at zero. We shall show that the condition obtained by
Basse and Pedersen~\cite{Ba} is also sufficient and give a totally
different proof for the necessity assertion, which is based on the
stationary increments property of fractional L\'evy processes. We
also obtain further characterisations based on differentiability
properties of $M_d$, show that the total variation is finite if and
only if its expectation is finite,  and obtain a zero-one law for
the property of being of finite variation. Note that when $M_d$ is a
semimartingale it may be used as a driving process for various
stochastic differential equations, and hence allows to incorporate
the long memory property into various classes of processes.

To set notation, fix a complete probability space $(\Omega,\cF,P)$
on which a real valued, two sided L\'evy process $L=(L(t))_{t\in
\bR}$ is defined, i.e. a process with independent and stationary
increments having c\`adl\`ag paths satisfying $L_0=0$. We shall
further assume throughout that $L(1)$ has finite variance and mean
zero. Recall that a two sided L\'evy processes $L$ indexed by $\bR$
can be easily constructed from a one-sided L\'evy processes $L_1$,
indexed by $[0,\infty)$, by letting $L(t)=L_1(t)$ for $t\geq 0$ and
$L(t) =  -L_2(-t-)$ for $t<0$, where $L_2$ is an independent copy of
$L_1$. We shall use the L\'evy--Khintchine representation of $L$ in
the form \begin{equation*} \label{eq-LK} \bE
e^{izL(1)}=\exp\left\{-\frac{1}{2}z^2\sigma+iz\gamma+\int_\bR
\left[e^{izx}-1-izx\One_{[-1,1]}(x)(1-|x|)\right]\,\nu(\mathrm
dx)\right\} \end{equation*} for  $z\in \mathbb{R}$, with
$\gamma\in\bR$, $\sigma\geq 0$ and $\nu$ being the L\'evy measure of
$L$, and refer to $(\sigma,\nu,\gamma)$ as the characteristic
triplet of $L$. See Sato~\cite{Sa} for further information regarding
L\'evy processes.

Let $d\in (0,1/2)$ which corresponds to a Hurst index $H:= d +1/2
\in (1/2,1)$ and hence to the long memory situation, to which we
limit ourselves in this paper. As in Marquardt~\cite{Ma}, define the
{\it non-anticipative fractional L\'evy process} $M_d =
(M_d(t))_{t\in \bR}$ by
\begin{eqnarray}
M_d(t)&:=&\frac{1}{\Gamma(d+1)}\int_\bR\big[(t-s)_+^d-(-s)_+^d\big]\
L(\mathrm ds), \quad  t\in\bR,\label{def-flp}
\end{eqnarray}
where for $\alpha\in \bR$ we put $x_+^\alpha := 0$ for $x\leq 0$ and
$x_+^\alpha := x^\alpha$ for $x > 0$. The integral in
\eqref{def-flp} converges in the $L^2$ sense (e.g.~\cite{Ma},
Proposition~2.1). As shown in~\cite{Ma}, Theorem~3.4, $M_d$ admits a
modification with continuous sample paths, which is given by the
following improper Riemann integral representation
\begin{equation}
\,\!M_d(t)=\frac{1}{\Gamma(d)}\int_{\bR}\!\left[(t-s)_+^{d-1}-(-s)_+^{d-1}\right]\!\,L(s)\,\mathrm
ds,\quad t\in\bR .\label{defMd}
\end{equation}
We shall always assume that we are working with the continuous
modification given by \eqref{defMd}.

Following  Samorodnitsky and Taqqu \cite{SaTa} the process defined
in (\ref{def-flp}) is called non-anticipative, since for positive
$t$ the value $M_d(t)$ depends on $\{L(s);s\leq t\}$ only. This is
in contrast to the {\it well-balanced fractional L\'evy process}
defined by
\begin{equation}
N_d(t):=\frac{1}{\Gamma(d+1)}\int_\bR\left(|t-s|^d-|s|^d\right)\
L(\mathrm ds),\quad  t\in\bR,\label{def-Nd}
\end{equation}
where $N_d(t)$ depends also on the future behaviour of $L$. It is
easy to see that both $M_d$ and $N_d$ have stationary increments,
which will be a crucial tool in the proofs presented.

The structure of the paper is as follows: Section~2 contains the
main results for the non-anticipative fractional L\'evy process
$M_d$, giving various characterisations for it to be of finite
variation,  one of which is the already mentioned condition in terms
of the characteristic triplet of $L$, while some others are given in
terms of differentiability properties of $M_d$.  A formula for the
derivative of $M_d$ at zero and for the expected total variation is
further obtained.  The results of Section~2 are
 proved in
Section~3. Section~4 treats the well-balanced fractional L\'evy
process $N_d$ and shows that the same characterisations hold for
$N_d$. Finally, in Section~5, the connection between fractional
L\'evy processes and the fractional Riemann-Liouville integral of a
L\'evy process is investigated, which establishes a direct link
between the results obtained in this article and those of Basse and
Pedersen~\cite{Ba} for this situation.

\section{Main results}
\setcounter{equation}{0}
The results in this section completely characterise when
fractional L\'evy processes are semimartingales. One characterisation is given
by integrability conditions on the L\'evy measure of the driving L\'evy process, a second by finiteness of the expected total variation, others are provided via various differentiability conditions on the sample paths of $M_d$. The characterisation also contains a 0-1 law for the finite variation property of the sample paths.

\begin{theorem} \label{thm-char}
Let $L$ be a L\'evy process with finite variance, $\bE L(1) = 0$ and
characteristic triplet $(\sigma,\nu, \gamma)$, let $d\in (0,1/2)$
and $[a,b]$ be a non-empty non-degenerate compact interval in $\bR$.
The following are equivalent:\\
(a) $M_d$ is almost surely of finite variation on $[a,b]$.\\
(a') $M_d$ is of finite variation on
$[a,b]$ with positive probability.\\
(b) With probability one the sample paths of $M_d$ are Lebesgue almost everywhere differentiable on $[a,b]$.\\
(b') With positive probability the sample paths of $M_d$ are
Lebesgue almost
everywhere differentiable on $[a,b]$.\\
(c) $M_d$ is almost surely differentiable at 0.\\
(c') $M_d$ is differentiable at 0 with positive probability.\\
(d) $M_d$ is almost surely differentiable from the left at 0.\\
(e) The Brownian motion part of $L$ is zero (i.e. $\sigma=0$) and
$$\int_{-1}^1 |x|^{\frac{1}{1-d}} \, \nu(\mathrm dx) < \infty.$$
(f) The expected total variation of $M_d$ on compacts is finite,
i.e.
$\bE\left[\mathrm{TV}\left(M_{d}|_{[a',b']}\right)\right]<\infty$
for all $a',b'\in \mathbb{R}$ with $a'<b'$.\\ (g) $M_d$ is a
semimartingale on $[0,\infty]$ with respect to every
filtration (satisfying the usual hypothesis) it is adapted to.\\
(h) $M_d$ is a semimartingale on $[0,\infty]$ with respect to some
filtration (satisfying the usual hypothesis) it is adapted to.
\end{theorem}
\begin{remark}
{\rm By Th\'eor\`eme III in Bretagnolle \cite{Br} property (e) in
the above theorem is equivalent to the driving L\'evy process $L$
being of finite $1/(1-d)$-variation.}
\end{remark}

\begin{remark} \label{rem-Basse-Pedersen}
{\rm Basse and Pedersen~\cite{Ba} showed that property (a) implies
(e) in the above theorem. This is stated explicitly in Corollary~5.4
of~\cite{Ba} when $\sigma$ is assumed a priori to be zero, but that
also $\sigma=0$ is necessary can be seen immediately from
Corollary~3.3 and Lemma~5.2 in \cite{Ba} along with a symmetrization
argument. Observe that in Theorem~\ref{thm-char} above, (e) is
deduced from the weaker conditions (b), (c), or (d), and that (e) is
shown to imply (a) and even the stronger condition (f).}
\end{remark}

For the derivative of $M_d$ at zero and for the expected total
variation, we have the following explicit representation:
\begin{theorem} \label{thm-derivative}
If the assumptions as well as one of the equivalent conditions of
Theorem~\ref{thm-char} hold, then the derivative of $M_d$ at zero can be expressed by
\begin{equation} \label{eq-derivative}
\frac{\mathrm d}{\mathrm dt} M_d(0) =\frac{1}{\Gamma(d)}
\int_{-\infty}^0 (-s)^{d-1} L(\mathrm ds) = \frac{1}{\Gamma(d-1)}
\int_{-\infty}^0 (-s)^{d-2} L(s) \, \mathrm ds,
\end{equation}
where  the integrals exist as improper integrals at zero in the sense of almost sure and $L^1(P)$ convergence. Moreover, for every $a<b\in \bR$,
$$
\bE\left[\mathrm{TV}\left(M_{d}|_{[a,b]}\right)\right]=\frac{(b-a)}{\Gamma(d)} \bE\left[\left|\int_{-\infty}^0 (-s)^{d-1} L(\mathrm ds)\right|\right].
$$
\end{theorem}

\section{The proofs}
\setcounter{equation}{0}

We first provide the proof of Theorem~\ref{thm-char} by showing the following implications
$$
(a)\Rightarrow (a') \Rightarrow (b') \Rightarrow (c') \Rightarrow
(c) \Rightarrow (d) \Rightarrow (e) \Rightarrow (f) \Rightarrow (g)
\Rightarrow (h) \Rightarrow (a)
$$
and
$$
(a) \Rightarrow (b) \Rightarrow (b').
$$
Note that
$$
(a)\Rightarrow (a') \Rightarrow (b'),\; (c) \Rightarrow (d),\; (f)
\Rightarrow (g) \Rightarrow (h),\; (a) \Rightarrow (b) \Rightarrow
(b')
$$
are obvious and $(h)\Rightarrow (a)$ is a consequence of the zero
quadratic variation property and the continuous paths of $M_d$, see
Theorem~4.7  in Marquardt \cite{Ma} for details.

\begin{proof}[Proof of Theorem~\ref{thm-char}, (b') $\Longrightarrow$ (c')]
By the stationary increments
property, we can assume without loss of generality that $[a,b] =
[-1,1]$. Denote
$$B := \{ (t_0, \omega) \in (-1,1) \times \Omega : \frac{\mathrm d}{\mathrm dt}
M_d(t_0,\omega) \; \mbox{exists}\}.$$ Note that the sample paths of
$M_d$ are continuous. Hence, $M_d:\bR\times \Omega\rightarrow \bR$
is $(\cB\otimes\cF)-\cB$-measurable,  and so is $h^{-1}
(M_d(t+h)-M_d(t))$ for every  $h\neq 0$, where $\cB$ denotes the
Borel $\sigma$-algebra in $\bR$. Thus, thanks to the continuous
paths of $M_d$, we observe that \begin{eqnarray*} B& =&
\left\{(t,\omega)\in (-1,1)\times\Omega:\
-\infty<\liminf_{h\in\mathbb{Q},\;h\rightarrow0}\frac{M_d(t+h)-M_d(t)}{h} \right.\\
& & \quad \quad \quad \quad \quad \quad \quad \quad \quad \quad
\left. =
\limsup_{h\in\mathbb{Q},\;h\rightarrow0}\frac{M_d(t+h)-M_d(t)}{h}
<\infty \right\}\in \cB\otimes\cF. \end{eqnarray*}
 Now, define the
cuts
$$B_{t_0} := \{\omega \in \Omega : (t_0,\omega) \in B \}, \quad
B^\omega := \{t_0 \in (-1,1) : (t_0,\omega) \in B\},$$ which are
$\cF$- and $\cB$-measurable, respectively. By assumption we have $P(
\{ \omega \in \Omega : \lambda_1 (B^\omega) = 1\})
> 0$, where $\lambda_1$ denotes the Lebesgue measure on $[-1,1]$. Further, by stationarity, $P(B_t) = P(B_0)$ for all $t \in
(-1,1)$, and consequently Fubini's theorem gives
$$2P(B_0) = \int_{-1}^1 P(B_t) \, \lambda_1(\mathrm dt) = (\lambda_1
\otimes P)\,\, (B) = \int_{\Omega} \lambda_1(B^\omega)\, dP > 0.$$ \end{proof}

\begin{proof}[Proof of Theorem~\ref{thm-char}, (c') $\Longrightarrow$ (c)]
With the definition above we have
$$B_0 = \{ \omega \in \Omega : \frac{\mathrm d}{\mathrm dt}
M_d(0,\omega) \; \mbox{exists}\}.$$ For fixed $r<0$ and $t\in (r,|r|)$, (\ref{defMd}) yields
$$\Gamma(d) M_d(t) =\! \int_{-\infty}^r \!\left[ (t-s)^{d-1} -
(-s)^{d-1}\right]L(s) \, \mathrm ds + \! \int_{r}^{0\vee t}\! \left[
(t-s)_+^{d-1} - (-s)_+^{d-1}\right]L(s) \, \mathrm ds . $$ Note that
$M_d(0) = 0$ and
$$\limsup_{s\to-\infty} \frac{|L(s)|}{(2 |s| \log \log |s|)^{1/2}} <
\infty$$ (cf. Sato~\cite{Sa}, Proposition~48.9). Moreover,  $L$ is
pathwise bounded on compacts. Therefore Lebesgue's dominated
convergence theorem ensures the existence of
$$\lim_{t\rightarrow 0} t^{-1} \int_{-\infty}^r \left[ (t-s)^{d-1}  - (-s)^{d-1}\right] L(s) \, \mathrm ds
= (d-1) \int_{-\infty}^r (-s)^{d-2} L(s) \,\mathrm ds.$$ It follows
that for each fixed $r < 0$, the set $B_0$ is measurable with
respect to $(L(s))_{r\leq s \leq |r|}$. Letting $r\uparrow 0$, we
conclude from the Blumenthal zero-one law, e.g. as given in
Proposition~40.4 in Sato~\cite{Sa}, that $P(B_0) \in \{0,1\}$. Since
$P(B_0) > 0$ by assumption, we have $P(B_0) = 1$.
\end{proof}

\begin{proof}[Proof of Theorem~\ref{thm-char}, (d) $\Longrightarrow$ (e)]
(i) Suppose first that $L$ is symmetric. Denote by $D$ the left
derivative of $\Gamma(d) M_d$ at 0. Then $D$ is infinitely divisible
as a limit of infinitely divisible distributions. Since $D$ is also
symmetric, its characteristic triplet with respect to the continuous
cut-off function $\beta(x) = (1 - |x|) \mathbf{1}_{[-1,1]}(x)$ is
given by $(A_D,\nu_D,0)$ with Gaussian variance $A_D$ and L\'evy
measure $\nu_D$, i.e. we have
$$\bE e^{iDz} = \exp \left[ -\frac12 A_D z^2 + \int_{\bR} \left(e^{izx}-1-i
zx \beta(x) \right) \, \nu_D(\mathrm dx)\right], \quad z\in \bR.$$
For any $r< 0$ we have (as shown in the proof of ``(c')
$\Longrightarrow$ (c)'') that
$$D = (d-1) \int_{-\infty}^r (-s)^{d-2} L(s) \, \mathrm ds + \lim_{t\uparrow
0} t^{-1} \int_{r}^0 \left[ (t-s)_+^{d-1} - (-s)_+^{d-1}\right] L(s) \,
\mathrm ds.$$
Integrating by parts, we obtain
\begin{eqnarray*}
(d-1) \int_{-\infty}^r (-s)^{d-2} L(s) \, \mathrm ds&=&  \int_{-\infty}^r (-s)^{d-1} L(\mathrm ds)-|r|^{d-1}L(r)\, , \\
t^{-1} \int_{r}^0 \left[ (t-s)_+^{d-1} - (-s)_+^{d-1}\right] L(s) \,
\mathrm ds &=& (td)^{-1}  \int_{r}^0 \left[ (t-s)_+^{d} -
(-s)_+^{d}\right] L(\mathrm ds) \\ &&+L(r)\frac{(|r|+t)^d-|r|^d}{td}
\, .
\end{eqnarray*}
Hence, for every $r<0$ we have
$$
D=Y_r+Z_r,
$$
where
\begin{eqnarray*}
Y_r & := & \int_{-\infty}^r (-s)^{d-1} \, L(\mathrm ds),\label{def-Y}\\
Z_r &:= & \lim_{t\uparrow 0} (td)^{-1} \int_{r}^0 \left[ (t-s)_+^{d}
- (-s)_+^{d}\right]\, L(\mathrm ds). \label{def-Z}
\end{eqnarray*}
Then $Y_r$ and $Z_r$ are independent for each $r$, and also
infinitely divisible. Denote the characteristic triplets of $Y_r$
and $Z_r$ with respect to $\beta$ by $(A_r,\nu_r,0)$ and $(A_{Z_r},
\nu_{Z_r},0)$, respectively, and observe that
\begin{equation*}
A_r+A_{Z_r}=A_D \mbox{ and } \nu_r + \nu_{Z_r} = \nu_D\label{A-nu-bounds}
\end{equation*}
by independence. Observe further that $A_r$ is a monotone increasing
sequence of real numbers bounded by $A_D$ and $(\nu_r)_{r<0}$ is an
increasing sequence of L\'evy measures as $r\uparrow 0$, bounded by
$\nu_D$. Denote
$$\nu_0(\Lambda) = \lim_{r\uparrow 0} \nu_r(\Lambda)$$
for each Borel set $\Lambda$. Then $\nu_0$ is a measure (e.g.
Kallenberg~\cite{Ka}, Corollary~1.16), and it is a L\'evy measure,
since it is bounded by $\nu_D$. Further,
$$\int_{[-\varepsilon,\varepsilon]} x^2 \nu_r(\mathrm dx) \leq
\int_{[-\varepsilon, \varepsilon]} x^2 \, \nu_D(\mathrm dx), \quad
\varepsilon > 0,\quad r < 0,$$ so that
$$\lim_{\varepsilon \downarrow 0} \limsup_{r\uparrow 0}
\int_{[-\varepsilon,\varepsilon]} x^2 \, \nu_r(\mathrm dx) = 0.$$ We
conclude that $Y_r$ converges in distribution to an infinitely
divisible random variable, $Y_0$ say, with characteristic triplet
$(A_0,\nu_0,0)$ with respect to $\beta$, cf. Sato~\cite{Sa},
Theorem~8.7. By the independent increments property of $r\mapsto
Y_r$, this convergence is even almost surely (e.g.
Kallenberg~\cite{Ka}, Theorem~4.18). Hence,
$$Y_{-1,0}:=\lim_{r\uparrow 0} \int_{-1}^r (-s)^{d-1} \, L(\mathrm ds)$$
exists as an almost sure limit. We now apply Proposition~5.3 parts (ii) and (i) of Sato~\cite{Sa07} to $f(s)=(-s)^{d-1}$ and $\alpha=1-\frac{d}{1-d}\in(0,1)$ to find that the existence of $\lim_{r\uparrow0}Y_{-1,r}$ in probability is equivalent to $L$ being purely non-Gaussian with a L\'evy measure $\nu$ fulfilling $\int_{[-1,1]} |x|^{1/(1-d)}\, \nu(\mathrm dx) < \infty$.\\
(ii) In the  general case denote by $\widetilde{L}$ an independent
copy of $L$, and write $L^* = L - \widetilde{L}$ for the
symmetrization of $L$. Suppose $M_d L$ (the fractional L\'evy
process driven by $L$) is almost surely left-differentiable at 0.
Then $M_d\widetilde{L}$ and $M_d L^*$ are almost surely
left-differentiable at 0 as well. Further, $L^*$ has a Gaussian part
if and only if $L$ has, and the L\'evy measure $\nu^*$ of $L^*$ is
given by $\nu^* (\Lambda) = \nu(\Lambda) + \nu(-\Lambda)$ for every
Borel set $\Lambda$. It follows from (i) that $\int_{-1}^1
|x|^{1/(1-d)} \nu^*(\mathrm dx) < \infty$, which is clearly
equivalent to $\int_{-1}^1 |x|^{1/(1-d)} \, \nu(\mathrm dx) <
\infty$.
\end{proof}

We now prepare the proof of the remaining implication
``(e)$\Rightarrow$(f)'' of Theorem~\ref{thm-char} and, at the same
time, the proof of Theorem \ref{thm-derivative} by the following
lemmas:

\begin{lemma}
\label{ozzy2} A symmetric and infinitely divisible random variable
$X$ without Gaussian part and with L\'evy measure $\nu$ fulfills the
following inequality for $\varepsilon>0$:
\begin{eqnarray}
\bE|X|&\leq& \varepsilon + \frac{4}{\varepsilon} \int_0^\varepsilon x \nu([x,\infty)) \mathrm dx + 2\int_\varepsilon^\infty\nu\left([x,\infty)\right)\,\mathrm dx. \nonumber
\end{eqnarray}
\begin{proof}
Writing $X=X^{(1)} + X^{(2)}$, where $X^{(1)}$ has characteristic
triplet $(0,\nu\big|_{(-\varepsilon,\varepsilon)},0)$ and $X^{(2)}$
has characteristic triplet
$(0,\nu\big|_{\bR\setminus(-\varepsilon,\varepsilon)},0)$ we find
$\mbox{for all }\varepsilon>0$ that $\bE|X|\leq
\bE|X^{(1)}|+\bE|X^{(2)}|.$ By Chebyshev's inequality, symmetry, the
fact that $\Var (X^{(1)}) = \int_{(-\varepsilon, \varepsilon)} x^2
\, \nu (\mathrm{d} x)$ (cf. Sato~\cite{Sa}, Example 25.12) and
integration by parts we get
\begin{eqnarray}
\!\!\!\bE|X^{(1)}|&=&\int_0^\infty P(|X^{(1)}|>t)\,\mathrm dt\leq \varepsilon +\int_\varepsilon^\infty\frac{\Var(X^{(1)})}{t^2}\,\mathrm dt\nonumber\\
&=&\varepsilon+\frac{2}{\varepsilon}\int_{(0,\varepsilon)}
x^2\,\nu(\mathrm
dx)=\varepsilon+\frac{4}{\varepsilon}\int_{0}^\varepsilon
x\,\nu([x,\infty))\mathrm dx-2 \varepsilon \nu([\epsilon,\infty)).
\label{eqn1}
\end{eqnarray}
Moreover, denoting by $X^{(3)}$ the compound Poisson distribution
with L\'evy measure $\nu\big|_{[\varepsilon,\infty)}$, we clearly
have
\begin{equation}
\bE|X^{(2)}|\leq 2\bE\,X^{(3)}=2\int_{[\varepsilon,\infty)}x\,\nu(\mathrm dx)=2\varepsilon\nu([\varepsilon,\infty))+2\int_{(\varepsilon,\infty)}\nu([x,\infty))\,\mathrm dx.\label{eqn2}
\end{equation}
Merging equations (\ref{eqn1}) and (\ref{eqn2}) gives the assertion.
\end{proof}
\end{lemma}

\begin{lemma}\label{vilbel}
Suppose condition (e) of Theorem \ref{thm-char} holds true. Then
$$\lim_{t\downarrow0}\bE\left(\left|\frac{1}{t}\int_0^t(t-s)^d\,L(\mathrm ds)\right|\right)=0.$$
\begin{proof}
(i) Let $L$ be symmetric and denote by $\nu_t$ the L\'evy measure of
$\frac{1}{t}\int_0^t(t-s)^d\,L(\mathrm ds).$ Then for every
$\varepsilon>0$ we obtain from Lemma~\ref{ozzy2} that
\begin{eqnarray}
\bE\left(\left|\frac{1}{t}\int_0^t(t-s)^d\,L(\mathrm ds)\right|\right)&\leq&\varepsilon+\frac{4}{\varepsilon}\int_0^\varepsilon u\nu_t([u,\infty))\,\mathrm du+2\int_\varepsilon^\infty\nu_t([u,\infty))\,\mathrm du\nonumber
\end{eqnarray}
Hence it is sufficient to show that, for every $\varepsilon>0$,
$$\mbox{a) } \lim_{t\downarrow0} \int_0^\varepsilon u\nu_t([u,\infty))\,\mathrm du=0,
\quad\mbox{and}\quad \mbox{b) } \lim_{t\downarrow0}
\int_\varepsilon^\infty\nu_t([u,\infty))\,\mathrm du=0.$$ To this
end note that, for $u>0$,
\begin{eqnarray*}\nu_t([u,\infty))&=&\int_0^t\int_{\bR\setminus\{0\}} \One_{[u,\infty)}\left(t^{-1}(t-s)^d x \right)\nu(\mathrm dx)\mathrm ds \\ &\leq& \int_0^t\int_{\{x>0\}} \One_{[u,\infty)}\left(t^{d-1} x \right)\nu(\mathrm dx)\mathrm ds=t\nu([ut^{1-d},\infty)),
\end{eqnarray*}
where the first identity is due to Theorem 3.10 in \cite{Sa07}.
 Hence,
\begin{equation}
0\leq \nu_t([u,\infty))\leq t\nu([ut^{1-d},\infty)) = \int_{ut^{1-d}}^\infty t\cdot x^\frac{1}{1-d}\cdot x^\frac{-1}{1-d}
\,\nu(\mathrm dx)\longrightarrow0\label{kl}
\end{equation}
by dominated convergence as $t\downarrow0$, since
$$
\One_{\{x\geq ut^{1-d}\}}\, t\cdot x^\frac{1}{1-d}\cdot x^\frac{-1}{1-d} \leq x^{\frac{1}{1-d}} u^{\frac{-1}{1-d}}.
$$
By (\ref{kl}) the integrand in a) converges to zero. So a) can be obtained by dominated convergence, because
$$u\nu_t([u,\infty))\leq u^{1-\frac{1}{1-d}}\int_0^\infty x^{\frac{1}{1-d}}\,\nu(\mathrm dx),\quad u>0,$$
and $u^{1-\frac{1}{1-d}}$ is integrable at zero since $d<\frac{1}{2}.$\\
Similarly, we get b) by (\ref{kl}) and dominated convergence, since
$$\int_{\varepsilon}^\infty \nu_t([u,\infty))\,\mathrm du\leq
t\int_{\varepsilon}^\infty\nu\left([ut^{1-d},\infty)\right)\,\mathrm
du \leq \int_0^\infty x^{\frac{1}{1-d}} \, \nu(\mathrm dx)
\int_\varepsilon^\infty u^{\frac{-1}{1-d}} \mathrm du < \infty.$$
(ii) If $L$ is not symmetric, then denote by $\tilde{L}$ an independent copy of $L$ and by $L^\ast=L-\tilde{L}$ the symmetrization of $L$. Further define
$$C_t=\frac{1}{t}\int_0^t(t-s)^d\,L(\mathrm ds)$$
and by $\tilde{C}_t,\,C^\ast_t$ analogous expressions for $\tilde L$ and $L^\ast$, respectively. As it holds
$$\bE\left[\left.\tilde{C}_t\right|C_t\right]=\bE\left[\tilde{C}_t\right]=0$$
by independence, we get
\begin{eqnarray}
\bE\left[|C_t|\right] &=& \bE\left[\left|C_t-\bE\left[\left.\tilde{C}_t\right|C_t\right]\right|\right]=\bE\left[\left|\bE\left[C_t^\ast|C_t\right]\right|\right]\leq \bE\left[|C_t^\ast|\right],\nonumber
\end{eqnarray}
and so the first case applies.
\end{proof}
\end{lemma}

\begin{lemma}
\label{ozzy1} Suppose condition (e) in Theorem \ref{thm-char} holds
true. For $r<0<t$ let $\nu_{r,t}$ denote the L\'evy measure of
$$B_{r,t}:=\frac{1}{t}\int_r^0\left[(t-s)^d-(-s)^d\right]\,L(\mathrm
ds).$$  Then for any $a> 0 $ it holds that
\begin{equation} \label{C2}
\int_0^au\cdot\nu_{r,t}\left([u,\infty)\right)\,\mathrm du \leq
\frac{a^2(1-d)}{1-2d}\left[ |r|
\nu\left(\left[\frac{a}{d}|r|^{1-d},\infty\right)\right) +
\left(\frac{d}{a}\right)^{\frac{1}{1-d}}
\int_0^{\frac{a}{d}|r|^{1-d}}x^{\frac{1}{1-d}}\,\nu(\mathrm
dx)\right]
\end{equation}
and
\begin{equation} \label{C3}
\int_a^\infty\nu_{r,t}\left([u,\infty)\right)\,\mathrm du\leq
|r|^d\int_{\frac{a}{d}|r|^{1-d}}^\infty x\, \nu(\mathrm
dx)+(1-d)\left(\frac{d}{a}\right)^{\frac{d}{1-d}}
\int_0^{\frac{a}{d}|r|^{1-d}} x^{\frac{1}{1-d}}\,\nu(\mathrm dx).
\end{equation}

\begin{proof}
From Theorem~3.10 in \cite{Sa07}, it follows that for any $u>0$,
\begin{eqnarray}
\nu_{r,t}\left([u,\infty)\right) & = & \int_0^{|r|}\int_0^\infty
\One_{[u,\infty)}\left(\frac{1}{t}\left[(t+s)^d-s^d\right]x\right)\,\nu(\mathrm
dx)\,\mathrm ds \nonumber\\
&\leq& \int_{0}^{|r|}\int_0^\infty\One_{[u,\infty)}(xs^{d-1}d)\,\nu(\mathrm dx)\,\mathrm ds\nonumber\\
&\leq&\int_{0}^{|r|}\nu\left(\left[\frac{u}{d} s^{1-d},\infty\right)\right)\,\mathrm ds\label{C4}\\
&=&\frac{1}{1-d}\cdot\left(\frac{d}{u}\right)^{\frac{1}{1-d}}\int_0^{\frac{u}{d}|r|^{1-d}}x^\frac{d}{1-d}\nu\left([x,\infty)\right)\,\mathrm
dx.\label{C1}
\end{eqnarray}
Hence, we obtain for any $a>0$ that
\begin{eqnarray*}
\int_0^au\cdot\nu_{r,t}\left([u,\infty)\right)\,\mathrm du
 &\leq&\frac{1}{1-d}d^{\frac{1}{1-d}}\int_0^a
 \int_0^{\frac{u}{d}|r|^{1-d}}u^{\frac{d}{d-1}}
 x^\frac{d}{1-d}\nu\left([x,\infty)\right)\,\mathrm dx\,\mathrm du\nonumber\\
&=&\frac{1}{1-d}d^{\frac{1}{1-d}}\int_0^{\frac{a}{d}|r|^{1-d}}
\int_{x|r|^{d-1}d}^au^{\frac{d}{d-1}}x^{\frac{d}{1-d}}\nu\left([x,\infty)\right)\,
\mathrm du\,\mathrm dx\nonumber\\
&\leq& \frac{1}{1-2d}d^\frac{1}{1-d} a^{1+\frac{d}{d-1}}
\int_0^{\frac{a}{d}|r|^{1-d}}x^{\frac{d}{1-d}}\nu\left([x,\infty)\right)\,\mathrm
dx,
\end{eqnarray*}
from which \eqref{C2} follows using integration by parts, observing
that the boundary term vanishes at 0 by \eqref{kl}. Moreover, by
\eqref{C1}, we have
\begin{eqnarray*}
\lefteqn{\int_a^\infty\nu_{r,t}\left([u,\infty)\right)\,\mathrm du} \nonumber \\
&\leq& \frac{d^\frac{1}{1-d}}{1-d} \int_a^\infty u^{-\frac{1}{1-d}}\int_0^{\frac{u}{d}|r|^{1-d}}
x^{\frac{d}{1-d}}\nu\left([x,\infty)\right)\,\mathrm dx\,\mathrm du\nonumber\\
&=& \frac{d^{\frac{1}{1-d}}}{1-d}\int_0^\infty\int_{\max\{a,xd|r|^{d-1}\}}^\infty
u^{\frac{1}{d-1}}\,\mathrm du\,\nu\left([x,\infty)\right)x^\frac{d}{1-d}\,\mathrm dx\nonumber\\
&=&  d^{\frac{d}{1-d}}\int_0^\infty \max\{a,xd|r|^{d-1}\}^{1-\frac{1}{1-d}}\,\nu
\left([x,\infty) \right) x^{\frac{d}{1-d}}\,\mathrm dx\nonumber\\
&=&  \left(\frac{d}{a}\right)^{\frac{d}{1-d}}
\int_0^{\frac{a}{d}|r|^{1-d}} x^\frac{d}{1-d}\nu\left([x,\infty)\right)\,
\mathrm dx+|r|^d \int_{\frac{a}{d}|r|^{1-d}}^\infty\nu\left([x,\infty)\right)\,\mathrm dx\nonumber\\
&=& (1-d)\left(\frac{d}{a}\right)^{\frac{d}{1-d}}
\int_0^{\frac{a}{d}|r|^{1-d}}x^{\frac{1}{1-d}}\,\nu(\mathrm dx)
+|r|^d\int^\infty_{\frac{a}{d}|r|^{1-d}}x\,\nu(\mathrm dx) \\
& & - a|r|\nu\left(\left[\frac{a}{d}|r|^{1-d},\infty\right)\right),
\end{eqnarray*}
giving \eqref{C3}.
\end{proof}
\end{lemma}

\begin{lemma}
\label{ozzy3}
Suppose condition (e) of  Theorem~\ref{thm-char} holds true. \\
{(a)} Then $B_{r,t}$, as defined in Lemma~\ref{ozzy1}, satisfies
$$\lim_{r\uparrow0}\sup_{t>0}E\left[\left|B_{r,t}\right|\right]=0.$$
{(b)} Moreover, the improper integral
$$\int_{-\infty}^0(-s)^{d-1}\,L(\mathrm ds):=\lim_{r\uparrow0}\int_{-\infty}^r(-s)^{d-1}\,L(\mathrm ds)$$
exists as an $L^1$-limit and as an almost sure limit.
\begin{proof}
(a) Note that, analogously to part (ii) in the proof of
Lemma~\ref{vilbel}, we may and do assume symmetry of $L$. For fixed
$\varepsilon>0$ we conclude from Lemmas~\ref{ozzy2} and~\ref{ozzy1}
that there is a constant $c_{d,\varepsilon}$ depending on $d$ and
$\varepsilon$ only, such that
\begin{eqnarray*}
\sup_{t>0} E[|B_{r,t}|]&\leq& \frac{\varepsilon}{2} + c_{d,\varepsilon}\left(|r|\nu\left([\varepsilon|r|^{1-d}/(2d),\infty)\right)+\int_0^{|r|^{1-d}\varepsilon/(2d)}x^{\frac{1}{1-d}}\nu(\mathrm dx)\right.\nonumber\\
&&\left.+|r|^d\int_{|r|^{1-d}\varepsilon/(2d)}^\infty x\nu(\mathrm dx)\right).
\end{eqnarray*}
The right-hand side is bounded by $\varepsilon$ for $r<0$
sufficiently close to 0, because
$$\lim_{r\uparrow0}\int_0^{|r|^{1-d}\varepsilon/(2d)}x^{1/(1-d)}\,\nu(\mathrm dx)=0$$
\begin{equation}
\lim_{r\uparrow0} |r|^d\int_{|r|^{1-d}\varepsilon/(2d)}x\nu(\mathrm dx)=0\label{second}
\end{equation}
\begin{equation}
\lim_{r\uparrow0}|r|\nu([|r|^{1-d},\infty))=0\label{third}
\end{equation}
by the assumption $\int_0^1x^{\frac{1}{1-d}}\,\nu(\mathrm dx)<\infty$ and dominated convergence. For (\ref{second}) the use of dominated convergence is justified by
\begin{eqnarray}
\One_{\{|r|^{1-d}\varepsilon/(2d)\leq x\leq1\}}|r|^dx&\leq&\One_{\{|r|^{1-d}\varepsilon/(2d)\leq x\leq1\}}|x|^\frac{1}{1-d}(\varepsilon/(2d))^{-d/(1-d)}\nonumber\\
&\leq&\One_{\{0<x\leq1\}}|x|^\frac{1}{1-d}(\varepsilon/(2d))^{-\frac{d}{1-d}},\nonumber
\end{eqnarray}
 Eq. (\ref{third}) is already shown in (\ref{kl}).\\
(b) The almost sure convergence follows from the $L^1$-convergence
by the independent increments property of $\int_{-\infty}^r
(-s)^{d-1}\,L(\mathrm ds),\,r<0,$ (see Theorem~4.18 in \cite{Ka}).
For proving the $L^1$-convergence, as above we may and do assume
that $L$ is symmetric.  Denote by $\tilde\nu_{q,r}$ the L\'evy
measure of
$$d\int_r^q(-s)^{d-1}\,L(\mathrm ds), \ \ \ \ \ \ -\infty<r<q<0.$$
Then, for $u>0$,
\begin{eqnarray}
\tilde\nu_{q,r}\left([u,\infty)\right)&=& \int_{|q|}^{|r|}\int_0^\infty\One_{[u,\infty)}(dxs^{d-1})\,\nu(\mathrm dx)\,\mathrm ds
\leq \int_0^{|r|} \nu\left([\frac{u}{d}s^{1-d},\infty)\right)\,\mathrm ds, \nonumber
\end{eqnarray}
which is the same upper bound as the one obtained for $\nu_{r,t}$ in
\eqref{C4} in the proof of Lemma \ref{ozzy1}. Therefore exactly the
same estimates as in part (a) can be applied to conclude that, given
$\varepsilon>0$, $ \bE\left[\left|\int_{r}^q d (-s)^{d-1}\,L(\mathrm
ds)\right|\right]\leq \varepsilon $ if $|r|, |q|$ are sufficiently
small. Hence  $\left(\int_{-\infty}^r (-s)^{d-1}\,L(\mathrm
ds)\right)_{r\uparrow0}$ is a Cauchy sequence and
therefore convergent in $L^1$. \\
\end{proof}
\end{lemma}

\begin{lemma}\label{ozzy4}
Under condition (e) of Theorem \ref{thm-char}, we have
$$
L^1-\lim_{t\downarrow 0} \frac{1}{t} M_d(t)=\frac{1}{\Gamma(d)} \int_{-\infty}^0(-s)^{d-1}\,L(\mathrm ds).
$$
\end{lemma}
\begin{proof}
 By Lemma~\ref{ozzy3} part b) the candidate limit is an $L^1$-random variable. Moreover, for every
 $r<0$ and $t>0$,
\begin{eqnarray}
&&\!\!\!\!\!\!\!\!\!\!\!\!\!\!\!\!\!\bE\left[\left|\frac{1}{t}\int_{-\infty}^t\left[(t-s)_+^d-(-s)^d_+\right]\,L(\mathrm ds) - d\int_{-\infty}^0(-s)^{d-1}\,L(\mathrm ds)\right|\right]\nonumber\\
&\leq& \bE\left[\left| \frac{1}{t}\int_{-\infty}^r\left[(t-s)^d-(-s)^d\right]\,L(\mathrm ds)-d\int_{-\infty}^r(-s)^{d-1}\,L(\mathrm ds)\right|\right]\nonumber\\
&+&\bE\left[\left|\frac{1}{t}\int_0^t(t-s)^d\,L(\mathrm ds)\right|\right]+\sup_{t'>0}\bE\left[\left|\frac{1}{t'}\int_r^0\left[(t'-s)^d-(-s)^d\right]\,L(\mathrm ds)\right|\right]\nonumber\\
&+&\bE\left[\left|d\int_r^0(-s)^{d-1}\,L(\mathrm ds)\right|\right]
=: (1) + (2) + (3) + (4), \quad \mbox{say.}\nonumber
\end{eqnarray}
Using the $L^2$-isometry for integrals with respect to square
integrable L\'evy processes (see e.g. Marquardt [7], Proposition
2.1), we get
\begin{eqnarray}
(1)&\leq& \left\{\bE\left\{\left|\int_{-\infty}^r\left\{\frac{1}{t}\left[(t-s)^d-(-s)^d\right]-d(-s)^{d-1}\right\}\,L(\mathrm ds)\right|^2\right\}\right\}^{1/2}\nonumber\\
&=&
\left(\bE\left(L(1)^2\right)\right)^{1/2}\left\{\int_{-\infty}^r\left\{\frac{1}{t}\left[(t-s)^d-(-s)^d\right]-d(-s)^{d-1}\right\}^2\mathrm
ds\right\}^{1/2}\nonumber
\end{eqnarray}
and by dominated convergence the latter integral converges to zero as $t\downarrow0$. For $t\downarrow0$ (2) tends to zero as shown in Lemma~\ref{vilbel}. Term (3) tends to zero as $r\uparrow 0$ by  Lemma~\ref{ozzy3}(a), and so does
term (4) by Lemma~\ref{ozzy3}(b). Hence, the assertion follows by letting $t\downarrow 0$ and then $r\uparrow 0$.
\end{proof}
\begin{remark}\label{rem_L1}
 Under condition (e) of Theorem \ref{thm-char}, one can also show that
$$
L^1-\lim_{t\uparrow 0} \frac{1}{t} M_d(t)=\frac{1}{\Gamma(d)} \int_{-\infty}^0(-s)^{d-1}\,L(\mathrm ds).
$$
To this end one decomposes, for $r<t<0$,
\begin{eqnarray*}
&&\!\!\!\!\!\!\!\!\!\!\!\!\!\!\!\!\!\bE\left[\left|\frac{1}{t}\int_{-\infty}^0\left[(t-s)_+^d-(-s)^d_+\right]\,L(\mathrm ds) - d\int_{-\infty}^0(-s)^{d-1}\,L(\mathrm ds)\right|\right]\nonumber\\
&\leq& \bE\left[\left| \frac{1}{t}\int_{-\infty}^r\left[(t-s)^d-(-s)^d\right]\,L(\mathrm ds)-d\int_{-\infty}^r(-s)^{d-1}\,L(\mathrm ds)\right|\right]\nonumber\\
&+&\bE\left[\left|\frac{1}{t}\int_t^0 (-s)^d\,L(\mathrm ds)\right|\right]+\sup_{r<t'<0}\bE\left[\left|\frac{1}{t'}\int_r^{t'}\left[(t'-s)^d-(-s)^d\right]\,L(\mathrm ds)\right|\right]\nonumber\\
&+&\bE\left[\left|d\int_r^0(-s)^{d-1}\,L(\mathrm ds)\right|\right]
\end{eqnarray*}
and shows convergence of these four terms analogously to the
situation in Lemma  \ref{ozzy4} letting $t\uparrow 0$ and then
$r\uparrow 0$.
\end{remark}

\begin{proof}[Proof of Theorem \ref{thm-char} `$(e)\Rightarrow (f)$', and Theorem \ref{thm-derivative}]
We assume that condition (e) in Theorem \ref{thm-char} is satisfied.
In order to prove (f) we fix $a'<b'\in\bR$ and for notational
convenience write $a=a'$ and $b=b'$. Due to the continuous paths of
$M_d$ the total variation can be calculated along dyadic partitions,
i.e.
$$
\mathrm{TV}\left(M_{d}|_{[a,b]}\right)=\lim_{n\rightarrow \infty} \sum_{i=1}^{2^n} |M_d(t_i)-M_d(t_{i-1})|,
$$
where $t_i=a+(b-a)i2^{-n}$. By the stationary increments of $M_d$, monotone convergence and Lemma \ref{ozzy4}, we obtain,
\begin{eqnarray*}
\bE\left[ \mathrm{TV}\left(M_{d}|_{[a,b]}\right)\right]&=&(b-a)\lim_{n\rightarrow \infty} \frac{2^n}{b-a}\bE[|M_d((b-a)2^{-n})|]\\&=&\frac{(b-a)}{\Gamma(d)} \bE\left[\left|\int_{-\infty}^0 (-s)^{d-1} L(\mathrm ds)\right|\right].
\end{eqnarray*}
In view of Lemma \ref{ozzy3}, the expectation on the right-hand side
is finite. Hence, (f) follows and the proof of Theorem
\ref{thm-char} is complete. Moreover, the explicit expression for
the expected total variation in Theorem \ref{thm-derivative} is
derived under condition (e) and, thus, under every of the equivalent
conditions in Theorem \ref{thm-char}.

Now suppose that one of the equivalent conditions in Theorem
\ref{thm-char} holds true. Then, properties (c) and (e) are valid.
By (c), the limit $M_d(t)/t$ exists almost surely as $t$ goes to
zero. Thanks to (e) and Lemma \ref{ozzy4} the limit $M_d(t)/t$
exists in $L^1$ as $t\downarrow 0$. Then both limits must coincide,
and consequently,
$$
\frac{\mathrm d}{\mathrm dt} M_d(0) =\frac{1}{\Gamma(d)}
\int_{-\infty}^0 (-s)^{d-1} L(\mathrm ds).
$$
Note that by Lemma \ref{ozzy3} the improper integral on the right hand side converges in $L^1$ and almost surely.
The alternative expression in (\ref{eq-derivative}) can be derived as follows. Integration by parts yields, for $r<0$,
$$
\int_{-\infty}^ r (-s)^{d-1} L(\mathrm ds)-(d-1)\int_{-\infty}^ r L(s)(-s)^{d-2}\mathrm ds=L(r)|r|^{d-1}.
$$
So, in view of Lemma \ref{ozzy3}, it is sufficient to show that the
right-hand side converges to zero in $L^1$ and almost surely for
$r\uparrow 0$. $L^1$-convergence is analogous to the proof of Lemma
\ref{vilbel} thanks to the integrability condition (e) for the
L\'evy measure, because the L\'evy measure $\tilde \nu_r$ of
$L(r)|r|^{d-1}$ is given by $\tilde \nu_r(\mathrm dx)=|r|
\nu(r^{1-d} \mathrm dx)$. Condition (e) also guarantees almost sure
convergence to zero in view of Theorem 2.1 in Bertoin, Doney and
Maller \cite{BDM08}.
\end{proof}

\section{Well-balanced fractional L\'evy processes}
\setcounter{equation}{0}

In this section we discuss the semimartingale property of the
closely related well-balanced fractional L\'evy process as defined
in \eqref{def-Nd}. Integration by parts yields the following
representation of a well-balanced fractional L\'evy process as
improper Riemann integral
\begin{equation}\label{Nbyparts}
 N_d(t)=\frac{1}{\Gamma(d)}\int_{-\infty}^{+\infty}\left(\sign(t-s)|t-s|^{d-1}+\sign(s)|s|^{d-1}\right)\,L(s)\mathrm
 ds .
\end{equation}
For well-balanced fractional L\'evy processes the equivalences of
Theorem~\ref{thm-char} do also hold true.

\begin{theorem} \label{thm-char2}
Under the assumption of Theorem \ref{thm-char} conditions (a)-(h)
given there are equivalent for the well-balanced fractional L\'evy
process $N_d$. If one of the conditions holds true, then
\begin{equation}
\frac{\mathrm d}{\mathrm dt} N_d(0) =\frac{-1}{\Gamma(d)}
\int_{-\infty}^\infty \sign(s)|s|^{d-1} L(\mathrm ds) = \frac{1}{\Gamma(d-1)}
\int_{-\infty}^\infty |s|^{d-2} L(s) \, \mathrm ds, \nonumber
\end{equation}
where the  integrals exist as improper integrals at zero in the
sense of almost sure and $L^1(P)$ convergence. Moreover, for every
$a'<b'\in \bR$,
$$
\bE\left[\mathrm{TV}\left(N_{d}|_{[a',b']}\right)\right]=\frac{(b'-a')}{\Gamma(d)}
\bE\left[\left|\int_{-\infty}^\infty \sign(s)|s|^{d-1} L(\mathrm
ds)\right|\right].
$$
\end{theorem}

The proof of Theorem \ref{thm-char2} follows similar lines as the proof of Theorem~\ref{thm-char}.
The following decomposition of well-balanced fractional L\'evy processes into the
 sum of two non-anticipative fractional L\'evy processes $M_d^{(1)}$ and $M_d^{(2)}$ turns out to be useful:
\begin{eqnarray}\label{NMdecomp}
&& \Gamma(d+1) N_d(t)= \int_{\bR}\left[(t-s)^d_+-(-s)_+^d+(t-s)_-^d-(-s)_-^d\right]\,L(\mathrm ds)\nonumber\\
& =& \int_{\bR} \left[(t-s)_+^d-(-s)_+^d\right]\,L(\mathrm ds) +
\int_{\bR}\left[(-t-u)^d_+-(-u)_+^d\right]\,L(-\mathrm du)\nonumber
\\ &=:& \Gamma(d+1)(M_d^{(1)}(t)+M_d^{(2)}(-t)).
\end{eqnarray}

\begin{proof}[Proof of Theorem \ref{thm-char2}]
The implications
$$
(a)\Rightarrow (a') \Rightarrow (b'),\; (c) \Rightarrow (d),\; (f)
\Rightarrow (g) \Rightarrow (h),\; (a) \Rightarrow (b) \Rightarrow
(b')
$$
are again obvious. \\
``(h) $\Longrightarrow$ (a)'': In view of (\ref{NMdecomp}) $N_d$ inherits the continuous paths and the zero quadratic variation of $M_d^{(i)}$, $i=1,2$ and, hence, is a semimartingale only if its paths are of bounded variation with probability one. \\
``(b') $\Longrightarrow$ (c')'': Analogously to Theorem~\ref{thm-char} using stationarity of $N_d$.\\
``(c') $\Longrightarrow$ (c)'': Denote
$$\tilde{B}_0=\{\omega\in\Omega:\frac{\mathrm d}{\mathrm dt}N_d(0,\omega) \mbox{ exists}\}.$$
Due to (\ref{Nbyparts}) we obtain, for $r>0$,
\begin{eqnarray*}
\Gamma(d)\cdot N_d(t)&=&\int_{-\infty}^{-r}\left(\sign(t-s)|t-s|^{d-1}+\sign(s)|s|^{d-1}\right)\,L(s)\mathrm ds\\&&+\int_{-r}^{r}\left(\sign(t-s)|t-s|^{d-1}+\sign(s)|s|^{d-1}\right)\,L(s)\mathrm ds\nonumber\\
&&+\int_{r}^{\infty}\left(\sign(t-s)|t-s|^{d-1}+\sign(s)|s|^{d-1}\right)\,L(s)\mathrm ds.\nonumber
\end{eqnarray*}
With the reasoning of the one-sided case it follows that the first and third limit in the corresponding decomposition of
$$\lim_{t\rightarrow0}\frac{1}{t}\Gamma(d)N_d(t)$$
exist for every fixed $r>0$. Hence $\tilde{B}_0$ is measurable with respect to $(L(s))_{-r\leq s \leq r}$, and the Blumenthal 0-1 law yields  $P(\tilde{B}_0)=1$.\\
``(d) $\Longrightarrow$ (e)'': Denoting the left derivative of
$\Gamma(d+1)N_d$ at $0$ by $\tilde{D}$ and using the representation
above we find for $r<t<0$ analogously to the non-anticipating case
\begin{eqnarray}
\tilde{D}&=&\int_{-\infty}^{r}|s|^{d-1}\,L(\mathrm ds)\nonumber \\ && + \lim_{t\uparrow0}(td)^{-1}\int_{r}^{|r|}\left(|t-s|^{d}-|s|^{d}\right)\,L(\mathrm ds)-\int_{|r|}^{+\infty}|s|^{d-1}\,L(\mathrm ds)\nonumber\\
&=:& Y_r^-+Z_r-Y_r^+.\nonumber
\end{eqnarray}
Then we apply the same reasoning to $Y_r^-$ as in the one-sided case for $Y_r$, giving almost sure convergence of $Y_r^-$ to a finite random variable as $r\uparrow 0$,
so that the claim follows again from Proposition~5.3 in Sato~\cite{Sa07}.\\
``(e) $\Longrightarrow$ (f)'': Fix $a'<b'\in\bR$. By the
decomposition (\ref{NMdecomp}) and the triangle inequality we have
$$
\bE\left[\mathrm{TV}\left(N_{d}|_{[a',b']}\right)\right]\leq
\bE\left[\mathrm{TV}\left(M^{(1)}_{d}|_{[a',b']}\right)\right]+\bE\left[\mathrm{TV}\left(M^{(2)}_{d}|_{[-b',-a']}\right)\right]<\infty
$$
thanks to Theorem \ref{thm-char}, (e) $\Longrightarrow$ (f).
\\[0.2cm]
Representation for the derivative and expected total variation: Applying Lemma \ref{ozzy4} and Remark \ref{rem_L1}  to $M_d^{(1)}$ and $M_d^{(2)}$ we obtain
\begin{eqnarray*}
 L^1-\lim_{t\rightarrow 0} \frac{1}{t} N_d(t)&=& L^1-\lim_{t\rightarrow 0} \frac{1}{t} M^{(1)}_d(t)+ L^1-\lim_{t\rightarrow 0} \frac{1}{t} M_d^{(2)}(-t)\\&=&\frac{1}{\Gamma(d)} \int_{-\infty}^0(-s)^{d-1}\,L(\mathrm ds)-\frac{1}{\Gamma(d)} \int_{-\infty}^0(-s)^{d-1}\,L(-\mathrm ds) \\&=& \frac{-1}{\Gamma(d)} \int_{-\infty}^\infty \sign(s)|s|^{d-1}\,L(\mathrm ds).
\end{eqnarray*}
Then an analogous reasoning as in the proof of Theorem \ref{thm-derivative} applies.
\end{proof}

\section{Fractionally integrated L\'evy processes}

Suppose $L$ is a L\'evy process with zero expectation and finite
variance. By (\ref{defMd}), the fractional L\'evy process $M_d$ of
order $0<d<1/2$ driven by $L$ can be split into the sum
$$
M_d(t)=\frac{1}{\Gamma(d)}\int_{-\infty}^0\!\left[(t-s)^{d-1}-(-s)^{d-1}\right]\!\,L(s) \mathrm ds +  (\cI^d L)(t) ,\quad t\geq 0.
$$
Here $(\cI^d f)(t)$ is the well-known fractional Riemann-Liouville integral of order $d>0$ defined by
$$\cI^d f(t)=\frac{1}{\Gamma(d)}\int_0^t(t-s)^{d-1}\,f(s)\,\mathrm ds$$
for sufficiently integrable functions $f$.

In order to transfer our results from fractional L\'evy processes to Riemann-Liouville integrals of L\'evy processes we first study the process
$$
F_d(t):=\frac{1}{\Gamma(d)}\int_{-\infty}^0\!\left[(t-s)^{d-1}-(-s)^{d-1}\right]\!\,L(s) \mathrm ds, \quad t\geq 0.
$$

\begin{proposition}
\label{carolo}
Suppose $L$ is a L\'evy process with zero expectation and finite variance and $0<d<1/2$. Then
 the expected total variation of $F_d$ on compact intervals $[0,b]$, $b>0$ is finite. Further,
 $\cI^d L$ is a.s. of finite variation on compact intervals (resp. has finite expected total variation on compact intervals) if and only if $M_d$ has the respective property.
\end{proposition}
\begin{proof}
Let $0=t_0 < t_1< \ldots < t_n=b$ be a partition of $[0,b]$. Then we
have for $s < 0 \leq t_{i-1} < t_i$ that $(t_i -s)^{d-1} <
(t_{i-1}-s)^{d-1}$ and we conclude
\begin{eqnarray*}
\lefteqn{\Gamma(d) \left| F_d(t_i) - F_d(t_{i-1})\right|}\\
& = & \left| \int_{-\infty}^0 \left[ (t_i-s)^{d-1} -
(t_{i-1}-s)^{d-1} \right] L(s) \, \mathrm ds\right| \\
& \leq &  \int_{-\infty}^0 \left[ (t_{i-1}-s)^{d-1} -
(t_i-s)^{d-1}\right] \, |L(s)| \, \mathrm ds.
\end{eqnarray*}
Taking in the following the supremum over all finite partitions of
$[0,b]$, the total variation
$\mathrm{TV}\left(F_{d}|_{[0,b]}\right)$ of $F_d$ over $[0,b]$ can
be estimated by
\begin{eqnarray*}
\mathrm{TV}\left(F_{d}|_{[0,b]}\right) & = & \sup_{t_0,\ldots, t_n}
\sum_{i=1}^n
|F_d(t_i) - F_d(t_{i-1})| \\
& \leq  & \frac{1}{\Gamma(d)} \sup_{t_0,\ldots, t_n}
\int_{-\infty}^0 \sum_{i=1}^n \left( (t_{i-1}-s)^{d-1} - (t_i
-s)^{d-1} \right) |L(s)| \, \mathrm ds \\
& = & \frac{1}{\Gamma(d)} \int_{-\infty}^0 \left[ (-s)^{d-1} -
(b-s)^{d-1} \right] |L(s)| \, \mathrm ds.
\end{eqnarray*}
Hence,
\begin{eqnarray*}
 \bE\left[\mathrm{TV}\left(F_{d}|_{[0,b]}\right)\right] &\leq& \frac{1}{\Gamma(d)} \int_{-\infty}^0 \left[ (-s)^{d-1} -
(b-s)^{d-1} \right] [E(|L(s)|^2)]^{1/2} \, \mathrm ds \\
&=& \frac{[E(|L(1)|^2)]^{1/2}}{\Gamma(d)} \int_{-\infty}^0 \left[
(-s)^{d-1} - (b-s)^{d-1} \right] s^{1/2} \, \mathrm ds
\end{eqnarray*}
The latter integral is finite because $0<d<1/2$ and
$$(-s)^{d-1} - (b-s)^{d-1} \sim b |d-1| |s|^{d-2} \quad
\mbox{as}\quad s\to -\infty.$$
Hence we see
that $F_d$ has finite expected total variation on compacts. The other assertions are immediate consequences of the decomposition
$$
M_d(t)=F_d(t)+(\cI^d L)(t),\quad t\geq 0.
$$
\end{proof}

Note that the fractional Riemann-Liouville integral $\cI^d L$ can be
defined for any L\'evy process $L$ (without any extra requirements
on the expectation and variance). In this general setting, the
combination of Corollaries 3.3 and 3.5(3) in Basse and Pedersen
\cite{Ba} states that, for $0<d<1/2$, $\cI^d L$ is a.s. of finite
variation on compacts, if and only if $L$ has no Gaussian component
and $|x|^{1/(1-d)}$ is integrable with respect to the L\'evy measure
$\nu$ around the origin (condition (e) in Theorem \ref{thm-char}
above). Hence, on the one hand, Proposition \ref{carolo} and the
results of Basse and Pedersen can be combined to provide an
alternative proof for the equivalence of properties (a) and (e) in
Theorem \ref{thm-char}. This observation complements Corollary 5.4
of \cite{Ba}, which only states that the condition (e) on the L\'evy
measure is necessary for the finite variation property of the paths
of the fractional L\'evy process $M_d$.

On the other hand, we can combine Theorem \ref{thm-char} and
Proposition \ref{carolo} to provide an alternative proof of Basse
and Pedersen's Corollary 3.5(3) and, additionally, include a 0-1 law
and a statement about the expected total variation, which is done in
the following result.

\begin{theorem}
Let $L$ be a L\'evy process with characteristic triplet
$(\sigma,\nu, \gamma)$, but without any moment assumptions.
Let $0<d<1/2$ and  $b>0$. Then the following statements are equivalent: \\
(a) $\cI^d L$ is a.s. of finite variation on $[0,b]$. \\
(b) $\cI^d L$ is of finite variation on $[0,b]$ with positive probability. \\
(c) The Brownian motion part of $L$ is zero (i.e. $\sigma = 0$) and
$$\int_{-1}^1 |x|^{\frac{1}{1-d}} \, \nu(\mathrm dx) < \infty.$$
Moreover, if one of the conditions holds true and $E[|L(1)|]<\infty$, then the expected total variation of $\cI^d L$ on $[0,b]$ is finite.
\end{theorem}
\begin{proof}
We decompose $L$ into a sum $L=L^{(1)}+L^{(2)}$, where the L\'evy
process $L^{(1)}$ contains the Gaussian part of $L$ and  the
compensated small jumps of $L$ corresponding to $\nu|_{(-1,1)}$.
Then $L^{(1)}$ has zero expectation and finite variance, and so
Theorem \ref{thm-char} and Proposition \ref{carolo} can be applied
to this process. Observe that
\begin{equation*}
\cI^d L=\cI^d L^{(1)}+\cI^d L^{(2)}=M_d^{(1)}-F_d^{(1)}+\cI^d
L^{(2)}.
\end{equation*}
In view of Theorem \ref{thm-char} and Proposition \ref{carolo} it remains to show that $\cI^d L^{(2)}$ is a.s. of finite variation on $[0,b]$, and that its expected total variation is finite provided $E[|L(1)|]<\infty$.  Note, that the L\'evy process $L^{(2)}$ only contains large jumps and a drift component and, thus, is of finite variation. We decompose $L^{(2)}=L^{(2,+)}-L^{(2,-)}$ into the difference of two increasing L\'evy processes (subordinators) $L^{(2,+)} $ and  $L^{(2,-)}$. Then,
$$
\cI^d L^{(2,\pm)}(t)=\frac{1}{\Gamma(d)}\int_0^t s^{d-1}\,L^{(2,\pm)}(t-s) \,\mathrm ds
$$
is increasing, because the mapping $t\mapsto {\bf 1}_{[0,t]}(s)s^{d-1}\,L^{(2,\pm)}(t-s)$ is increasing for fixed $s$. Then, we obtain, for every $b>0$,
$$
\mathrm{TV}\left(\cI^d L^{(2)}|_{[0,b]}\right)\leq  \cI^d
L^{(2,+)}(b)+\cI^d L^{(2,-)}(b)<\infty.
$$
Note that the L\'evy measures of $L^{(2,\pm)}$ are  $\nu|_{[1,\infty)}$ and $\nu|_{(-\infty,-1]}$. Hence, $E[|L(1)|]<\infty$ implies that $E[|L^{(2,+)}(t)|+|L^{(2,-)}(t)|]<\infty$ for every $t>0$. Thus,
$$
E\left[\mathrm{TV}\left(\cI^d L^{(2)}|_{[0,b]}\right)\right]\leq  \frac{1}{\Gamma(d)}\int_0^b s^{d-1}\,E[|L^{(2,+)}(b)|+|L^{(2,-)}(b)|]\,\mathrm ds<\infty.
$$
\end{proof}

\subsection*{Acknowledgements}

We are indebted to Muneya Matsui for valuable comments and fruitful
discussions, which in particular started this project. Thanks also
to Ren\'e Schilling for his help with the proof of the implication
``(b') $\Longrightarrow$ (c')'' of Theorem~2.1. Further thanks go to
the referee for careful reading and valuable suggestions. Support
from an NTH-grant of the state of Lower Saxony is gratefully
acknowledged.

\end{document}